# On competing risk and degradation processes

**Nozer D. Singpurwalla**[1],[*]

*The George Washington University*

**Abstract:** Lehmann's ideas on concepts of dependence have had a profound effect on mathematical theory of reliability. The aim of this paper is two-fold. The first is to show how the notion of a "hazard potential" can provide an explanation for the cause of dependence between life-times. The second is to propose a general framework under which two currently discussed issues in reliability and in survival analysis involving interdependent stochastic processes, can be meaningfully addressed via the notion of a *hazard potential*. The first issue pertains to the failure of an item in a dynamic setting under multiple interdependent risks. The second pertains to assessing an item's life length in the presence of observable surrogates or markers. Here again the setting is dynamic and the role of the marker is akin to that of a leading indicator in multiple time series.

## 1. Preamble: Impact of Lehmann's work on reliability

Erich Lehmann's work on non-parametrics has had a conceptual impact on reliability and life-testing. Here two commonly encountered themes, one of which bears his name, encapsulate the essence of the impact. These are: the notion of a **Lehmann Alternative**, and his exposition on **Concepts of Dependence**. The former (see Lehmann [4]) comes into play in the context of accelerated life testing, wherein a Lehmann alternative is essentially a model for accelerating failure. The latter (see Lehmann [5]) has spawned a large body of literature pertaining to the reliability of complex systems with interdependent component lifetimes. Lehmann's original ideas on characterizing the nature of dependence has helped us better articulate the effect of failures that are *causal* or *cascading*, and the consequences of lifetimes that exhibit a negative correlation. The aim of this paper is to propose a framework that has been inspired by (though not directly related to) Lehmann's work on dependence. The point of view that we adopt here is "dynamic", in the sense that what is of relevance are dependent stochastic processes. We focus on two scenarios, one pertaining to competing risks, a topic of interest in survival analysis, and the other pertaining to degradation and its markers, a topic of interest to those working in reliability. To set the stage for our development we start with an overview of the notion of a *hazard potential*, an entity which helps us better conceptualize the process of failure and the cause of interdependent lifetimes.

---

[*]Research supported by Grant DAAD 19-02-01-0195, The U. S. Army Research Office.
[1]Department of Statistics, The George Washington University, Washington, DC 20052, USA, e-mail: nozer@gwu.edu

## 2. Introduction: The hazard potential

Let $T$ denote the time to failure of a unit that is scheduled to operate in some specified static environment. Let $h(t)$ be the hazard rate function of the survival function of $T$, namely, $P(T \geq t), t \geq 0$. Let $H(t) = \int_0^t h(u)du$, be the cumulative hazard function at $t$; $H(t)$ is increasing in $t$. With $h(t), t \geq 0$ specified, it is well known that

$$\Pr(T \geq t; h(t), t \geq 0) = \exp(-H(t)).$$

Consider now an exponentially distributed random variable $X$, with scale parameter $\lambda$, $\lambda \geq 0$. Then for some $H(t) \geq 0$,

$$\Pr(X \geq H(t) | \lambda = 1) = \exp(-H(t));$$

thus

(2.1) $\qquad \Pr(T \geq t; h(t), t \geq 0) = \exp(-H(t)) = \Pr(X \geq H(t) | \lambda = 1).$

The right hand side of the above equation says that the item in question will fail when its cumulative hazard $H(t)$ crosses a threshold $X$, where $X$ has a unit exponential distribution. Singpurwalla [11] calls $X$ the **Hazard Potential** of the item, and interprets it as an unknown resource that the item is endowed with at inception. Furthermore, $H(t)$ is interpreted as the amount of resource consumed at time $t$, and $h(t)$ is the rate at which that resource gets consumed. Looking at the failure process in terms of an endowed and a consumed resource enables us to characterize an environment as being *normal* when $H(t) = t$, and as being *accelerated* (*decelerated*) when $H(t) \geq (\leq) t$. More importantly, with $X$ interpreted as an unknown resource, we are able to interpret dependent lifetimes as the consequence of dependent hazard potentials, the later being a manifestation of commonalities of design, manufacture, or genetic make-up. Thus one way to generate dependent lifetimes, say $T_1$ and $T_2$ is to start with a bivariate distribution $(X_1, X_2)$ whose marginal distributions are exponential with scale parameter one, and which is not the product of exponential marginals. The details are in Singpurwalla [11].

When the environment is dynamic, the rate at which an item's resource gets consumed is random. Thus $h(t); t \geq 0$ is better described as a stochastic process, and consequently, so is $H(t), t \geq 0$. Since $H(t)$ is increasing in $t$, the **cumulative hazard process** $\{H(t); t \geq 0\}$ is a continuous increasing process, and the item fails when this process hits a random threshold $X$, the item's hazard potential. Candidate stochastic processes for $\{H(t); t \geq 0\}$ are proposed in the reference given above, and the nature of the resulting lifetimes described therein. Noteworthy are an increasing Lévy process, and the maxima of a Wiener process.

In what follows we show how the notion of a hazard potential serves as a unifying platform for describing the competing risk phenomenon and the phenomenon of failure due to ageing or degradation in the presence of a marker (or a bio marker) such as crack size (or a CD4 cell count).

## 3. Dependent competing risks and competing risk processes

By "competing risks" one generally means failure due to agents that presumably compete with each other for an item's lifetime. The traditional model that has been used for describing the competing risk phenomenon has been the reliability of a series system whose component lifetimes are independent or dependent. The idea here is that since the failure of any component of the system leads to the failure of the system, the system experiences multiple risks, each risk leading to failure. Thus if $T_i$ denotes the lifetime of component $i$, $i = 1, \ldots, k$, say, then the cause of system failure is that component whose lifetime is smallest of the $k$ lifetimes. Consequently, if $T$ denotes a system's lifetime, then

$$\text{(3.1)} \qquad \Pr(T \geq t) = P(H_1(t) \leq X_1, \ldots, H_k(t) \leq X_k),$$

where $X_i$ is the hazard potential of the $i$-th component, and $H_i(t)$ its cumulative hazard (or the risk to component $i$) at time $t$. If the $X_i$'s are assumed to be independent (a simplifying assumption), then (3.1) leads to the result that

$$\text{(3.2)} \qquad \Pr(T \geq t) = \exp[-(H_1(t) + \cdots + H_k(t))],$$

suggesting an additivity of cumulative hazard functions, or equivalently, an additivity of the risks. Were the $X_i$'s assumed dependent, then the nature of their dependence will dictate the manner in which the risks combine. Thus for example if for some $\theta$, $0 \leq \theta \leq 1$, we suppose that

$$\Pr(X_1 \geq x_1, X_2 \geq x_2 | \theta) = \exp(-x_1 - x_2 - \theta x_1 x_2),$$

namely one of Gumbel's bivariate exponential distributions, then

$$\Pr(T \geq t | \theta) = \exp[-(H_1(t) + H_2(t) + \theta H_1(t) H_2(t))].$$

The cumulative hazards (or equivalently, the risks) are no longer additive.

The series system model discussed above has also been used to describe the failure of a single item that experiences several failure causing agents that compete with each other. However, we question this line of reasoning because a single item posseses only one unknown resource. Thus the $X_1, \ldots, X_k$ of the series system model should be replaced by a single $X$, where $X_1 = X_2 = \cdots = X_k = X$ (in probability). To set the stage for the single item case, suppose that the item experiences $k$ agents, say $C_1, \ldots, C_k$, where an agent is seen as a cause of failure; for example, the consumption of fatty foods. Let $H_i(t)$ be the consequence of agent $C_i$, were $C_i$ be the only agent acting on the item. Then under the simultaneous action by all of the $k$ agents the item's survival function

$$\text{(3.3)} \qquad \begin{aligned} \Pr(T \geq t; h_1(t), \ldots, h_k(t)) \\ = P(H_1(t) \leq X, \ldots, H_k(t) \leq X) \\ = \exp(-\max(H_1(t), \ldots, H_k(t))). \end{aligned}$$

Here again, the cumulative hazards are not additive.

Taking a clue from the fact that dependent hazard potentials lead us to a non-additivity of the cumulative hazard functions, we observe that the condition $X_1 \stackrel{P}{=} X_2 \stackrel{P}{=} \cdots \stackrel{P}{=} X_k \stackrel{P}{=} X$ (where $X_1 \stackrel{P}{=} X_2$ denotes that $X_1$ and $X_2$ are equal in probability) implies that $X_1, \ldots, X_k$ are *totally positively dependent*, in the sense of Lehmann (1966). Thus (3.2) and (3.3) can be combined to claim that in general, under the series system model for competing risks, $P(T \geq t)$ can be bounded as

$$\text{(3.4)} \qquad \exp(-\sum_1^k H_i(t)) \leq P(T \geq t) \leq \exp(-\max(H_1(t), \ldots, H_k(t))).$$

Whereas (3.4) above may be known, our argument leading up to it could be new.

### 3.1. Competing risk processes

The prevailing view of what constitutes dependent competing risks entails a consideration of dependent component lifetimes in the series system model mentioned above. By contrast, our position on a proper framework for describing dependent competing risks is different. Since it is the $H_i(t)$'s that encapsulate the notion of risk, dependent competing risks should entail interdependence between $H_i(t)$'s, $i = 1, \ldots, k$. This would require that the $H_i(t)$'s be random, and a way to do so is to assume that each $\{H_i(t); t \geq 0\}$ is a stochastic process; we call this a **competing risk process**. The item fails when any one of the $\{H_i(t); t \geq 0\}$ processes first hits the items hazard potential $X$. To incorporate interdependence between the $H_i(t)$'s, we conceptualize a $k$-variate process $\{H_1(t), \ldots, H_k(t); t \geq 0\}$, that we call a **dependent competing risk process**. Since $H_i(t)$'s are increasing in $t$, one possible choice for each $\{H_i(t); t \geq 0\}$ could be a *Brownian Maximum Process*. That is $H_i(t) = \sup_{0 < s \leq t}\{W_i(s); s \geq 0\}$, where $\{W_i(s); s \geq 0\}$ is a standard Brownian motion process. Dependence between the $H_i(t)$'s can be induced via a dependence between the $\{W_i(s); s \geq 0\}$ processes. Thus for example, in the bivariate case, if $\rho$ denotes the correlation between two standard Brownian motion processes, then

$$\Pr(T \geq t) = \int_0^\infty P\left(H_1(t) \leq x, H_2(t) \leq x\right) e^{-x} dx$$

and it can be shown (details omitted) that,

$$(3.5) \quad \Pr(T \geq t) = \frac{\int_0^t \int_0^t \exp\left[-\frac{(a^2+b^2-2\rho ab)}{2t(1-\rho^2)}\right] da\, db}{\int_0^\infty \int_0^\infty \exp\left[-\frac{(u^2+v^2-2\rho uv)}{2t(1-\rho^2)}\right] du\, dv}.$$

Another possibility, again for the case of $k = 2$, is to assume that $\{H_1(t); t \geq 0\}$ is some non-negative, non-decreasing, right-continuous process, but that $\{H_2(t); t \geq 0\}$ has a sample path which is an impulse function of the form $H_2(t) = 0$ for all $t < t^*$, and that $H_2(t^*) = \infty$ for some $t^* > 0$, where the rate of occurrence of the impulse at time $t$ depends on $H_1(t)$. The process $\{H_2(t); t \geq 0\}$ can be identified with some sort of a traumatic event that competes with the process $\{H_1(t); t \geq 0\}$ for the lifetime of the item. In the absence of trauma the item fails when the process $\{H_1(t); t \geq 0\}$ hits the item's hazard potential. This scenario parallels the one considered by Lemoine and Wenocur [6], albeit in a context that is different from ours. By assuming that the probability of occurrence of an impulse in the time interval $[t, t+h)$, given that $H_1(t) = \omega$, is $1 - \exp(-\omega h)$, Lemoine and Wenocur [6] have shown that for $X = x$, the probability of survival of an item to time $t$ is of the form:

$$(3.6) \quad \Pr(T \geq t) = E\left[\exp\left(\int_0^t H_1(s) ds\right) I_{[0,x)}\left(H_1(t)\right)\right],$$

where $I_A(\bullet)$ is the indicator of a set $A$, and the expectation is with respect to the distribution of the process $\{H_1(t); t \geq 0\}$. As a special case, when $\{H_1(t); t \geq 0\}$ is a gamma process (see Singpurwalla [10]), and $x$ is infinite, so that $I_{[0,\infty)}\left(H_1(t)\right) = 1$ for $H_1(t) \geq 0$, the above equation takes the form

$$(3.7) \quad \Pr(T \geq t) = \exp(-(1+t)\log(1+t) + t).$$

The closed form result of (3.7) suffers from the disadvantage of having the effect of the hazard potential de facto nullified. The more realistic case of (3.6) will call for numerical or simulation based approaches. These remain to be done; our aim here has been to give some flavor of the possibilities.

## 4. Biomarkers and degradation processes

A topic of current interest in both reliability and survival analysis pertains to assessing lifetimes based on observable surrogates, such as crack length, and biomarkers like CD4 cell counts. Here again the hazard potential provides a unified perspective for looking at the interplay between the unobservable failure causing phenomenon, and an observable surrogate. It is an assumed dependence between the above two processes that makes this interplay possible.

To engineers (cf. Bogdanoff and Kozin [1]) degradation is the irreversible accumulation of damage throughout life that leads to failure. The term "damage" is not defined; however it is claimed that damage manifests itself via surrogates such as cracks, corrosion, measured wear, etc. Similarly, in the biosciences, the notion of "ageing" pertains to a unit's position in a state space wherein the probabilities of failure are greater than in a former position. Ageing manifests itself in terms of biomedical and physical difficulties experienced by individuals and other such biomarkers.

With the above as background, our proposal here is to conceptualize ageing and degradation as unobservable constructs (or latent variables) that serve to describe a process that results in failure. These constructs can be seen as the cause of observable surrogates like cracks, corrosion, and biomarkers such as CD4 cell counts. This modelling viewpoint is not in keeping with the work on degradation modelling by Doksum [3] and the several references therein. The prevailing view is that degradation is an observable phenomenon that reveals itself in the guise of crack length and CD4 cell counts. The item fails when the observable phenomenon hits some threshold whose nature is not specified. Whereas this may be meaningful in some cases, a more general view is to separate the observable and the unobservable and to attribute failure as a consequence of the behavior of the unobservable.

To mathematically describe the cause and effect phenomenon of degradation (or ageing) and the observables that it spawns, we view the (unobservable) cumulative hazard function as degradation, or ageing, and the biomarker as an observable process that is influenced by the former. The item fails when the cumulative hazard function hits the item's hazard potential $X$, where $X$ has exponential (1) distribution. With the above in mind we introduce the **degradation process** as a bivariate stochastic process $\{H(t), Z(t), t \geq 0\}$, with $H(t)$ representing the unobservable degradation, and $Z(t)$ an observable marker. Whereas $H(t)$ is required to be non-decreasing, there is no such requirement on $Z(t)$. For the marker to be useful as a predictor of failure, it is necessary that $H(t)$ and $Z(t)$ be related to each other. One way to achieve this linkage is via a *Markov Additive Process* (cf. Cinlar [2]) wherein $\{Z(t); t \geq 0\}$ is a Markov process and $\{H(t); t \geq 0\}$ is an increasing Lévy process whose parameters depend on the state of the $\{Z(t); t \geq 0\}$ process. The ramifications of this set-up need to be explored.

Another possibility, and one that we are able to develop here in some detail (see Section 5), is to describe $\{Z(t); t \geq 0\}$ by a Wiener process (cracks do heal and CD4 cell counts do fluctuate), and the unobservable degradation process $\{H(t); t \geq 0\}$ by a ***Wiener Maximum Process***, namely,

$$H(t) = \sup_{0 < s \leq t} \{Z(s); s \geq 0\}. \tag{4.1}$$

What makes the topic of analyzing degradation processes attractive is not just the modeling part; the statistical and computational issues that the set-up creates are quite challenging. Since $\{Z(t); t \geq 0\}$ is an observable process, how may one use observations on this process until some time, say $t^*$, to make inferences about the process of interest $H(t)$, for $t > t^*$? In other words, how does one assess $\Pr(T > t|\{Z(s); 0 < s \leq t^* < t\})$, where $T$ is an item's time to failure? Furthermore, as is often the case, the process $\{Z(s); s \geq 0\}$ cannot be monitored continuously. Rather, what one is able to do is observe $\{Z(s); s \geq 0\}$ at $k$ discrete time points and use these as a basis for inference about $\Pr(T > t|\{Z(s); 0 < s \leq t^* < t\})$. These and other matters are discussed next in Section 5, which could be viewed as a prototype of what else is possible using other models for degradation.

## 5. Inference under a Wiener maximum process for degradation

We start with some preliminaries about a Wiener process and its hitting time to a threshold. The notation used here is adopted from Doksum [3].

### *5.1. Hitting time of a Wiener maximum process to a random threshold*

Let $Z_t$ denote an observable marker process $\{Z(t); t \geq 0\}$, and $H_t$ an unobservable degradation process $\{H(t); t \geq 0\}$. The relationship between these two processes is prescribed by (4.1). Suppose that $Z_t$ is described by a Wiener process with a drift parameter $\eta$ and a diffusion parameter $\sigma^2 > 0$. That is, $Z(0) = 0$ and $Z_t$ has independent increments. Also, for any $t > 0$, $Z(t)$ has a Gaussian distribution with $E(Z(t)) = \eta t$, and for any $0 \leq t_1 < t_2$, $\text{Var}[Z(t_2) - Z(t_1)] = (t_2 - t_1)\sigma^2$. Let $T_x$ denote the <u>first</u> time at which $Z_t$ crosses a threshold $x > 0$; that is, $T_x$ is the hitting time of $Z_t$ to $x$. Then, when $\eta = 0$,

$$\begin{aligned}\Pr(Z(t) \geq x) &= \Pr(Z(t) \geq x | T_x \leq t)\Pr(T_x \leq t) \\ &\quad + \Pr(Z(t) \geq x | T_x > t)\Pr(T_x > t),\end{aligned} \tag{5.1}$$

$$\Pr(T_x \leq t) = 2\Pr(Z(t) \geq x). \tag{5.2}$$

This is because $\Pr(Z(t) \geq x | T_x \leq t)$ can be set to $1/2$, and the second term on the right hand side of (5.1) is zero. When $Z(t)$ has a Gaussian distribution with mean $\eta t$ and variance $\sigma^2 t$, $\Pr(Z(t) \geq x)$ can be similarly obtained, and thence $\Pr(T_x \leq t) \stackrel{def}{=} F_x(t|\eta, \sigma)$. Specifically it can be seen that

$$F_x(t|\eta, \sigma) = \Phi\left(\frac{\sqrt{\lambda}}{\mu}\sqrt{t} - \frac{\sqrt{\lambda}}{\sqrt{t}}\right) + \Phi\left(-\frac{\sqrt{\lambda}}{\mu}\sqrt{t} - \frac{\sqrt{\lambda}}{\sqrt{t}}\right)\exp\left(\frac{2\lambda}{\mu}\right), \tag{5.3}$$

where $\mu = x/\eta$ and $\lambda = x^2/\sigma^2$. The distribution $F_x$ is the ***Inverse Gaussian Distribution*** (*IG*-Distribution) with parameters $\mu$ and $\lambda$, where $\mu = E(T_x)$ and $\lambda \mu^2 = \text{Var}(T_x)$. Observe that when $\eta = 0$, both $E(T_x)$ and $\text{Var}(T_x)$ are infinite, and thus for any meaningful description of a marker process via a Wiener process, the drift parameter $\eta$ needs to be greater than zero.

The probability density of $F_x$ at $t$ takes the form:

(5.4) $$f_x(t|\eta,\sigma) = \sqrt{\frac{\lambda}{2\pi t^3}} \exp\left[-\frac{\lambda}{2\mu^2}\frac{(t-\mu)^2}{t}\right],$$

for $t, \mu, \lambda > 0$.

We now turn attention to $H_t$, the process of interest. We first note that because of (4.1), $H(0) = 0$, and $H(t)$ is non-decreasing in $t$; this is what was required of $H_t$. An item experiencing the process $H_t$ fails when $H_t$ first crosses a threshold $X$, where $X$ is unknown. However, our uncertainty about $X$ is described by an exponential distribution with probability density $f(x) = e^{-x}$. Let $T$ denote the time to failure of the item in question. Then, following the line of reasoning leading to (5.1), we would have, in the case of $\eta = 0$,

$$\Pr(T \leq t) = 2\Pr(H(t) \geq x).$$

Furthermore, because of (4.1), the hitting time of $H_t$ to a random threshold $X$ will coincide with $T_x$, the hitting time of $Z_t$ (with $\eta > 0$) to $X$. Consequently,

$$\Pr(T \leq t) = \Pr(T_x \leq t) = \int_0^\infty \Pr(T_x \leq t|X=x)f(x)dx$$
$$= \int_0^\infty \Pr(T_x \leq t)e^{-x}dx = \int_0^\infty F_x(t|\eta,\sigma)e^{-x}dx.$$

Rewriting $F_x(t|\eta,\sigma)$ in terms of the marker process parameters $\eta$ and $\sigma$, and treating these parameters as known, we have

(5.5) $$\Pr(T \leq t|\eta,\sigma) \stackrel{def}{=} F(t|\eta,\sigma)$$
$$= \int_0^\infty \left[\Phi\left(\frac{\eta}{\sigma}\sqrt{t} - \frac{x}{\sigma\sqrt{t}}\right) + \Phi\left(-\frac{\eta}{\sigma}\sqrt{t} - \frac{x}{\sigma\sqrt{t}}\right)\right]$$
$$\times \exp\left(x\left(\frac{2\eta}{\sigma^2} - 1\right)\right)dx,$$

as our assessment of an item's time to failure with $\eta$ and $\sigma$ assumed known. It is convenient to summarize the above development as follows

**Theorem 5.1.** *The time to failure $T$ of an item experiencing failure due to ageing or degradation described by a **Wiener Maximum Process** with a drift parameter $\eta > 0$, and a diffusion parameter $\sigma^2 > 0$, has the distribution function $F(t|\eta,\sigma)$ which is a location mixture of Inverse Gaussian Distributions. This distribution function, which is also the hitting time of the process to an exponential (1) random threshold, is given by (5.5).*

In Figure 1 we illustrate the behavior of the *IG*-Distribution function $F_x(t)$, for $x = 1, 2, 3, 4$, and 5, when $\eta = \sigma = 1$, and superimpose on these a plot of $F(t|\eta = \sigma = 1)$ to show the effect of averaging the threshold $x$. As can be expected, averaging makes the *S*-shapedness of the distribution functions less pronounced.

### *5.2. Assessing lifetimes using surrogate (biomarker) data*

The material leading up to Theorem 5.1 is based on the thesis that $\eta$ and $\sigma^2$ are known. In actuality, they are of course unknown. Thus, besides the hazard potential $X$, the $\eta$ and $\sigma^2$ constitute the unknowns in our set-up. To assess $\eta$ and $\sigma^2$ we may use prior information, and when available, data on the underlying processes $Z_t$ and $H_t$. The prior on $X$ is an exponential distribution with scale one, and this prior can also be updated using process data. In the remainder of this section, we focus attention on the case of a single item and describe the nature of the data that can be collected on it. We then outline an overall plan for incorporating these data into our analyses.

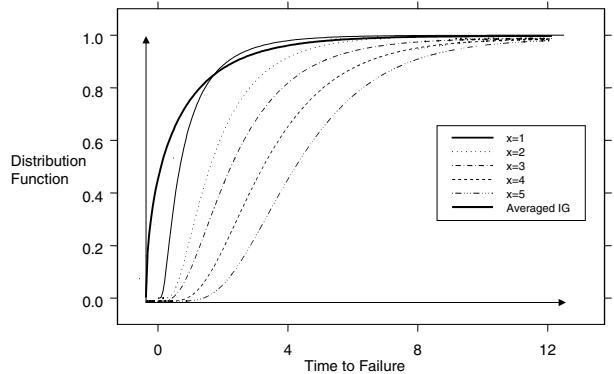

Fig 1. *The IG-Distribution with thresholds $x = 1, \ldots, 5$ and the averaged IG-Distribution.*

In Section 5.3 we give details about the inferential steps. The scenario of observing several items to failure in order to predict the lifetime of a future item will not be discussed.

In principle, we have assumed that $H_t$ is an unobservable process. This is certainly true in our particular case when the observable marker process $Z_t$ cannot be continuously monitored. Thus it is not possible to collect data on $H_t$. Contrast our scenario to that of Doksum [3], Lu and Meeker [7], and Lu, Meeker and Escobar [8], who assume that degradation is an observable process and who use data on degradation to predict an item's lifetime. We assume that it is the surrogate (or the biomarker) process $Z_t$ that is observable, but only prior to $T$, the item's failure time. In some cases we may be able to observe $Z_t$ at $t=T$, but doing so in the case of a single item would be futile, since our aim is to assess an unobserved $T$. Data on $Z_t$ will certainly provide information about $\eta$ and $\sigma^2$, but also about $X$; this is because for any $t < T$, we know that $X > Z(t)$. Thus, as claimed by Nair [9], data on (the observable surrogates of) degradation helps sharpen lifetime assessments, because a knowledge of $\eta$, $\sigma^2$ and $X$ translates to a knowledge of $T$.

It is often the case – at least we assume so – that $Z_t$ cannot be continuously monitored, so that observations on $Z_t$ could be had only at times $0 < t_1 < t_2 < \cdots < t_k < T$, yielding $\mathbf{Z} = (Z(t_1), \ldots, Z(t_k))$ as data. Furthermore, based on $Z(t_k)$, we are able to assert that $X > Z(t_k)$. This means that our updated uncertainty about $X$ will be encapsulated by a shifted exponential distribution with scale parameter one, and a location (or shift) parameter $Z(t_k)$.

Thus for an item experiencing failure due to degradation, whose marker process yields $\mathbf{Z}$ as data, our aim will be to assess the item's *residual life* $(T - t_k)$. That is, for any $u > 0$, we need to know $\Pr(T > t_k + u; \mathbf{Z}) = \Pr(T > t_k + u; T > t_k)$, and this under a certain assumption (cf. Singpurwalla [12]) is tantamount to knowing

$$(5.6) \qquad \frac{\Pr(T > t_k + u)}{\Pr(T > t_k)},$$

for $0 < u < \infty$. To assess the two quantities in the above ratio, we need to consider the quantity $\Pr(T > t; \mathbf{Z})$, for some $t > 0$. Let $\pi(\eta, \sigma^2, x; \mathbf{Z})$ encapsulate our uncertainty about $\eta$, $\sigma^2$ and $X$ in the light of the data $\mathbf{Z}$. In Section 5.3 we describe our approach for assessing $\pi(\eta, \sigma^2, x; \mathbf{Z})$. Now

$$\text{(5.7)} \quad \Pr(T > t; \mathbf{Z}) = \int_{\eta, \sigma^2, x} \Pr(T > t | \eta, \sigma^2, x; \mathbf{Z}) \pi(\eta, \sigma^2, x; \mathbf{Z})(d\eta)(d\sigma^2)(dx)$$

$$= \int_{\eta, \sigma^2, x} \Pr(T_x > t | \eta, \sigma^2) \pi(\eta, \sigma^2, x; \mathbf{Z})(d\eta)(d\sigma^2)(dx)$$

$$\text{(5.8)} \quad = \int_{\eta, \sigma^2, x} F_x(t | \eta, \sigma) \pi(\eta, \sigma^2, x; \mathbf{Z})(d\eta)(d\sigma^2)(dx),$$

where $F_x(t|\eta, \sigma)$ is the *IG*-Distribution of (5.3).

Implicit to going from (5.7) to (5.8) is the assumption that the event $(T > t)$ is independent of $\mathbf{Z}$ given $\eta$, $\sigma^2$ and $X$. In Section 5.3 we will propose that $\eta$ be allowed to vary between $a$ and $b$; also, $\sigma^2 > 0$, and having observed $Z(t_k)$, it is clear that $x$ must be greater than $Z(t_k)$. Consequently, (5.8) gets written as

$$\text{(5.9)} \quad \Pr(T > t; \mathbf{Z}) = \int_a^b \int_0^\infty \int_{Z(t_k)}^\infty F_x(t|\eta, \sigma) \pi(\eta, \sigma^2, x; \mathbf{Z})(d\eta)(d\sigma^2)(dx),$$

and the above can be used to obtain $\Pr(T > t_k + u; \mathbf{Z})$ and $\Pr(T > t_k; \mathbf{Z})$. Once these are obtained, we are able to assess the residual life $\Pr(T > t_k + u | T > t_k)$, for $u > 0$.

We now turn our attention to describing a Bayesian approach specifying $\pi(\eta, \sigma^2, x; \mathbf{Z})$.

### 5.3. Assessing the posterior distribution of $\eta$, $\sigma^2$ and $X$

The purpose of this section is to describe an approach for assessing $\pi(\eta, \sigma^2, x; \mathbf{Z})$, the posterior distribution of the unknowns in our set-up. For this, we start by supposing that $\mathbf{Z}$ is an unknown and consider the quantity $\pi(\eta, \sigma^2, x | \mathbf{Z})$. This is done to legitimize the ensuing simplifications. By the multiplication rule, and using obvious notation

$$\pi(\eta, \sigma^2, x|\mathbf{Z}) = \pi_1(\eta, \sigma^2|X, \mathbf{Z}) \pi_2(X|\mathbf{Z}).$$

It makes sense to suppose that $\eta$ and $\sigma^2$ do not depend on $X$; thus

$$\text{(5.10)} \quad \pi(\eta, \sigma^2, x|\mathbf{Z}) = \pi_1(\eta, \sigma^2|\mathbf{Z}) \pi_2(X|\mathbf{Z}).$$

However, $\mathbf{Z}$ is an observed quantity. Thus (5.10) needs to be recast as:

$$\text{(5.11)} \quad \pi(\eta, \sigma^2, x; \mathbf{Z}) = \pi_1(\eta, \sigma^2; \mathbf{Z}) \pi_2(X; \mathbf{Z}).$$

Regarding the quantity $\pi_2(X; \mathbf{Z})$, the only information that $\mathbf{Z}$ provides about $X$ is that $X > Z(t_k)$. Thus $\pi_2(X; \mathbf{Z})$ becomes $\pi_2(X; Z(t_k))$. We may now invoke Bayes' law on $\pi_2(X; Z(t_k))$ and using the facts that the prior on $X$ is an exponential (1) distribution on $(0, \infty)$, obtain the result that the posterior of $X$ is also an exponential (1) distribution, but on $(Z(t_k), \infty)$. That is, $\pi_2(X; Z(t_k))$ is a shifted exponential distribution of the form $\exp(-(x - Z(t_k)))$, for $x > Z(t_k)$.

Turning attention to the quantity $\pi_1(\eta, \sigma^2; \mathbf{Z})$ we note, invoking Bayes' law, that

$$\text{(5.12)} \quad \pi_1(\eta, \sigma^2; \mathbf{Z}) \propto \mathcal{L}(\eta, \sigma^2; \mathbf{Z}) \pi^*(\eta, \sigma^2),$$

where $\mathcal{L}(\eta, \sigma^2; \mathbf{Z})$ is the likelihood of $\eta$ and $\sigma^2$ with $\mathbf{Z}$ fixed, and $\pi^*(\eta, \sigma^2)$ our prior on $\eta$ and $\sigma^2$. In what follows we discuss the nature of the likelihood and the prior.

**The Likelihood of $\eta$ and $\sigma^2$**

Let $Y_1 = Z(t_1)$, $Y_2 = (Z(t_2) - Z(t_1)), \ldots, Y_k = (Z(t_k) - Z(t_{k-1}))$, and $s_1 = t_1$, $s_2 = t_2 - t_1, \ldots, s_k = t_k - t_{k-1}$. Because the Wiener process has independent increments, the $y_i$'s are independent. Also, $y_i \sim N(\eta s_i, \sigma^2 s_i)$, $i = 1, \ldots, k$, where $N(\mu, \xi^2)$ denotes a Gaussian distribution with mean $\mu$ and variance $\xi^2$. Thus, the joint density of the $y_i$'s, $i = 1, \ldots, k$, which is useful for writing out a likelihood of $\eta$ and $\sigma^2$, will be of the form

$$\prod_{i=1}^{k} \phi\left(\frac{y_i - \eta s_i}{\sigma^2 s_i}\right);$$

where $\phi$ denotes a standard Gaussian probability density function. As a consequence of the above, the likelihood of $\eta$ and $\sigma^2$ with $\mathbf{y} = (y_1, \ldots, y_k)$ fixed, can be written as:

$$(5.13) \qquad \mathcal{L}(\eta, \sigma^2; \mathbf{y}) = \prod_{i=1}^{k} \frac{1}{\sqrt{2\pi s_i}\sigma} \exp\left[-\frac{1}{2}\left(\frac{y_i - \eta s_i}{\sigma^2 s_i}\right)^2\right].$$

**The Prior on $\eta$ and $\sigma^2$**

Turning attention to $\pi^*(\eta, \sigma^2)$, the prior on $\eta$ and $\sigma^2$, it seems reasonable to suppose that $\eta$ and $\sigma^2$ are not independent. It makes sense to suppose that the fluctuations of $Z_t$ depend on the trend $\eta$. The larger the $\eta$, the bigger the $\sigma^2$, so long as there is a constraint on the value of $\eta$. If $\eta$ is not constrained the marker will take negative values. Thus, we need to consider, in obvious notation

$$(5.14) \qquad \pi^*(\eta, \sigma^2) = \pi^*(\sigma^2|\eta)\pi^*(\eta).$$

Since $\eta$ can take values in $(0, \infty)$, and since $\eta = \tan\theta$ – see Figure 2 – $\theta$ must take values in $(0, \pi/2)$.

To impose a constraint on $\eta$, we may suppose that $\theta$ has a *translated beta* density on $(a, b)$, where $0 < a < b < \pi/2$. That is, $\theta = a + (b-a)W$, where $W$ has a beta distribution on $(0, 1)$. For example, $a$ could be $\pi/8$ and $b$ could be $3\pi/8$. Note that were $\theta$ assumed to be uniform over $(0, \pi/2)$, then $\eta$ will have a density of the form $2/[\pi(1 + \eta^2)]$ – which is a *folded Cauchy*.

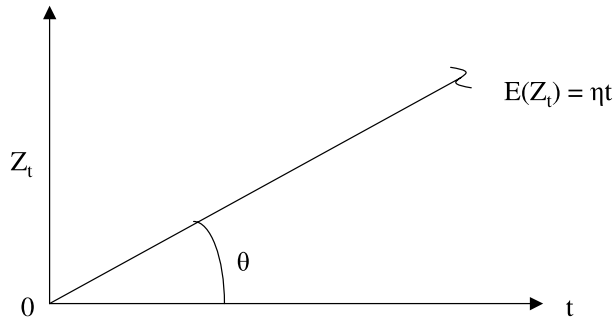

FIG 2. *Relationship between $Z_t$ and $\eta$.*

The choice of $\pi^*(\sigma^2|\eta)$ is trickier. The usual approach in such situations is to opt for natural conjugacy. Accordingly, we suppose that $\psi \stackrel{def}{=} \sigma^2$ has the prior

$$\pi^*(\psi|\eta) \propto \psi^{-\left(\frac{\nu}{2}+1\right)} \exp\left(-\frac{\eta}{2\psi}\right), \quad (5.15)$$

where $\nu$ is a parameter of the prior.

Note that $E(\psi|\eta,\nu) = \eta/(\nu-2)$, and so $\psi = \sigma^2$ increases with $\eta$, and $\eta$ is constrained over $a$ and $b$. Thus a constraint on $\sigma^2$ as well.

To pin down the parameter $\nu$, we anchor on time $t = 1$, and note that since $E(Z_1) = \eta$ and $\text{Var}(Z_1) = \sigma^2 = \psi$, $\sigma$ should be such that $\Delta\sigma$ should not exceed $\eta$ for some $\Delta = 1, 2, 3, \ldots$; otherwise $Z_1$ will become negative. With $\Delta = 3$, $\eta = 3\sigma$ and so $\psi = \sigma^2 = \eta^2/9$. Thus $\nu$ should be such that $E(\sigma^2|\eta,\nu) \approx \eta^2/9$. But $E(\sigma^2|\eta,\nu) = \eta/(\nu-2)$, and therefore by setting $\eta/(\nu-2) = \eta^2/9$, we would have $\nu = 9/\eta + 2$. In general, were we to set $\eta = \Delta\sigma$, $\nu = \Delta^2/\eta + 2$, for $\Delta = 1, 2, \ldots$. Consequently, $\nu/2 + 1 = (\Delta^2/\eta + 2)/2 + 1 = \Delta^2/2\eta + 2$, and thus

$$\pi^*(\psi|\eta;\Delta) = \psi^{-\left(\frac{\Delta^2}{2\eta}+2\right)} \exp\left(-\frac{\eta}{2\psi}\right), \quad (5.16)$$

would be our prior of $\sigma^2$, conditioned on $\psi$, and $\Delta = 1, 2, \ldots$, serving as a prior parameter. Values of $\Delta$ can be used to explore sensitivity to the prior.

This completes our discussion on choosing priors for the parameters of a Wiener process model for $Z_t$. All the necessary ingredients for implementing (5.9) are now at hand. This will have to be done numerically; it does not appear to pose major obstacles. We are currently working on this matter using both simulated and real data.

## 6. Conclusion

Our aim here was to describe how Lehmann's original ideas on (positive) dependence framed in the context of non-parametrics have been germane to reliability and survival analysis, and even so in the context of survival dynamics. The notion of a hazard potential has been the "hook" via which we can attribute the cause of dependence, and also to develop a framework for an appreciation of competing risks and degradation. The hazard potential provides a platform through which the above can be discussed in a unified manner. Our platform pertains to the hitting times of stochastic processes to a random threshold. With degradation modeling, the unobservable cumulative hazard function is seen as the metric of degradation (as opposed to an observable, like crack growth) and when modeling competing risks, the cumulative hazard is interpreted as a risk. Our goal here was not to solve any definitive problem with real data; rather, it was to propose a way of looking at two commonly encountered problems in reliability and survival analysis, problems that have been well discussed, but which have not as yet been recognized as having a common framework. The material of Section 5 is purely illustrative; it shows what is possible when one has access to real data. We are currently persuing the details underlying the several avenues and possibilities that have been outlined here.

## Acknowledgements

The author acknowledges the input of Josh Landon regarding the hitting time of a Brownian maximum process, and Bijit Roy in connection with the material of Section 5. The idea of using Wiener Maximum Processes for the cumulative hazard was the result of a conversation with Tom Kurtz.